\documentclass{amsart}
\usepackage{amssymb, amsmath, amscd, latexsym, mathrsfs, appendix}
\usepackage{graphics}
\usepackage{amsmath}
\usepackage{amsxtra}
\usepackage{amsfonts}
\usepackage{appendix}
\usepackage[all]{xy}

\usepackage{color}

\numberwithin{equation}{section}
\newtheorem{thm}{Theorem}[section]

\newtheorem{lem}[thm]{Lemma}

\newtheorem{defn}[thm]{Definition}

\newtheorem{expl}[thm]{Example}

\newcommand{\lra}{\longrightarrow}

\newcommand{\co}{\colon\!}

\newcommand{\smin}{\smallsetminus}

\newcommand{\id}{\textup{id}}

\newcommand{\holim}{\textup{holim}}

\newcommand{\hofiber}{\textup{hofiber}}
\newcommand{\mor}{\textup{mor}}

\newcommand{\rmap}{\mathbb R\textup{map}}

\newcommand{\config}{\mathsf{con}} 

\newcommand{\fin}{\mathsf{Fin}}
\newcommand{\dend}{\mathsf{Tree}}
\newcommand{\dendroca}{\dend^{\mathsf{rc}}}
\newcommand{\dendlev}{\mathsf{Levtree}}

\newcommand{\simp}{\mathsf{simp}}

\newcommand{\sC}{\mathcal C}

\newcommand{\sS}{\mathcal S}

\newcommand{\op}{\textup{op}}

\newcommand{\ms}{\mathcal}

\newcommand{\RR}{\mathbb R}

\newcommand{\holimsub}[1]{\begin{array}[t]{cc} \textup{holim} \\ [-1mm]
\scriptstyle{#1} \end{array}}

\newcommand{\uli}{\underline}

\begin{document}

\title{Truncated operads and simplicial spaces}
\author{Michael S. Weiss}%

\address{Math.Institut, WWU M\"{u}nster, 48149 M\"{u}nster, Einsteinstrasse 62, Germany}%
 \email{m.weiss@uni-muenster.de}

\thanks{The work was supported by the Bundesministerium f\"ur Bildung und Forschung through the A.v.Humboldt foundation
(Humboldt professorship 2012-2017).}

\begin{abstract} It was shown in \cite[\S7]{BoavidaWeissLong} that a well-known construction
which to a (monochromatic, symmetric) topological operad associates a topological category and a functor from it to the category of finite sets
is homotopically fully faithful, under mild conditions on the operads. The main result here is a
generalization of that statement to $k$-truncated topological operads. A $k$-truncated operad is a weaker version of
operad where all operations have arity $\le k$.
\end{abstract}
\maketitle

\section{Introduction}
There is a well known construction which to a (monochromatic, symmetric) topological operad $P$ associates a so-called PROP, a small topological category with a
symmetric monoidal product. The set of objects of the PROP is identified with the set of natural numbers
and the monoidal product corresponds to addition. See \cite{Voronov} for more details and explanations. Forgetting the monoidal product in
the PROP and passing to the comma category of \emph{objects over 1} produces a topological category (category object in the category of spaces)
called $\sC_P$ in \cite[\S7]{BoavidaWeissLong}. This comes with a distinguished functor to $\fin$, essentially the category of all finite sets.
It is well known that $P$ can be reconstructed from the associated PROP with the monoidal
structure. By contrast, the construction $P\mapsto \sC_P$ has a forgetful character, even if we include the reference functor $\sC_P\to \fin$,
and there are elementary examples to illustrate that it really does forget essential features \cite[rem.7.3]{BoavidaWeissLong}.
Write $N$ for the nerve construction (from small topological
categories to simplicial spaces). The main result of \cite[\S7]{BoavidaWeissLong} is that
for topological operads $P$ and $Q$, the map
\[  \rmap(P,Q) \lra \rmap_{N\fin}(N\sC_P,N\sC_Q) \]
of derived mapping spaces induced by the construction $P\mapsto N\sC_P$ is nevertheless a weak equivalence under a reasonable condition on $P$ and $Q$.
The condition is that the spaces of $0$-ary and $1$-ary operations for both $P$ and $Q$ are weakly contractible.
(Beware that topological operads with an \emph{empty space} of $0$-ary operations do not qualify, although they are popular.) More details and an example are given in section~\ref{sec-oper}.

\smallskip
This paper here relies on \cite[\S7]{BoavidaWeissLong}, but in doing so develops a slightly different proof
of the same result. The advantage of the new proof is that it carries over without essential change to the setting of
\emph{$k$-truncated} topological operads. (A $k$-truncated operad is a weaker version of
operad which has operations of arity $\le k$ only.)

\section{Operads, dendroidal spaces and simplicial spaces}  \label{sec-oper}
This section is a review of the main definitions and results of \cite[\S7]{BoavidaWeissLong}.

\smallskip
Let $P$ be an operad in the symmetric monoidal category of spaces. For convenience,
the \emph{category of spaces} is understood to be the category of simplicial sets. The following paragraph in quotation marks is
quoted verbatim from \cite[\S7]{BoavidaWeissLong}. \newline
``We think of this in the following terms: $P$ is a functor from the category of finite sets and bijections to
spaces, and for every map $f:T\to S$ of finite sets there is an operation
\[  \lambda_f\co P(S)\times\prod_{i\in S} P(T_i) \lra P(T) \]
where $T_i=f^{-1}(i)\subset T$. Also $P(S)$ contains a distinguished unit element when $S$ is a singleton.
Sensible naturality, associativity and unital properties are satisfied. Note in particular that any
permutation $f\co S\to S$ induces a map $P(S)\to P(S)$ in two ways: firstly because
$P$ is a functor from the category of finite sets and bijections to spaces, and secondly by
\[ P(S)\ni x\mapsto \lambda_f(x,1,1,\dots,1)\in P(S)~. \]\
for $x\in P(S)$. These two maps agree as per definition. \newline
What we have described is also called a \emph{plain} operad in the category of spaces ...''

Note that \emph{plain} is synonymous with \emph{monochromatic}. It means that the operad has only one object.
(We may think of $P(S)$ as the space of $S$-ary operations from that object to itself.) The fact that we
allow arbitrary finite sets without a total ordering means that
we are dealing with a \emph{symmetric} operad. In the following we just write \emph{operad} to mean
\emph{monochromatic symmetric operad}.

\smallskip
There is a construction $P\mapsto \sC_P$ taking an operad $P$ as above to a topological category $\sC_P$
(category object in the category of spaces). The category $\sC_P$ comes with a forgetful functor to the category
of finite sets; more precisely we use a skeleton $\fin$ of the category of finite sets. The objects of $\fin$ are the
sets $\uli k=\{1,2,\dots,k\}$ where $k\ge 0$. (Note in particular that $\uli 0$ is the empty set.) The morphisms in $\fin$
from $\uli k$ to $\uli\ell$ are the maps from $\uli k$ to $\uli\ell$. \newline
The space of objects of $\sC_P$ is
\[  \coprod_{k\ge 0} P(\uli k)\,. \]
The space of morphisms in $\ms C_P$ lifting a morphism $f\co\uli k\to \uli\ell$ in $\fin$ is
\[ P(\uli\ell)\times\prod_{i\in \uli\ell} P(\uli k_{\,i}) \]
where $k_i$ is the cardinality of $f^{-1}(i)$.
Source and target of an element in that space are determined by applying to it $\lambda_f$
and the projection to $P(\ell)$, respectively. Composition and identity morphisms are obvious. See \cite[Rem.7.3]{BoavidaWeissLong}
for the relationship between the construction $P\mapsto\sC_P$ and a better known construction \cite{Voronov} which turns an operad into a PROP.

\smallskip
The main result of \cite[\S7]{BoavidaWeissLong} states that the construction $P\mapsto \sC_P$ is homotopically fully faithful
as long as it is only used with operads $P$ for which the spaces $P(\uli 0)$ and $P(\uli 1)$ are weakly contractible. Making this precise
requires a few decisions. Clearly $P\mapsto \sC_P$ is a functor. But we need to reason about derived spaces of morphisms between two topological operads,
or between two small topological categories. Therefore it is convenient to work with a preferred model category of topological operads, and with a
preferred model category of small topological categories.

\smallskip
Rezk has introduced a model category setting for the \emph{category of small topological categories}.
The underlying category is simply the category of simplicial spaces (and here \emph{space} still means \emph{simplicial set} for us).
A simplicial space $X$ is said to be a \emph{Segal space} if it satisfies the following condition $(\sigma)$.
Let $u_i\co\{0,1\}\to\{0,1,2,\dots,n\}$ be the order-preserving map defined by $u_i(0)=i-1$ and $u_i(1)=i$.
\begin{itemize}
\item[($\sigma$)] For each $n\ge 2$ the map $(u_1^*,u_2^*,\dots,u_n^*)$ from $X_n$ to the homotopy limit of
the diagram
\[
\xymatrix{
X_1\ar[r]^-{d_0} & X_0 & \ar[l]_{d_1} X_1 \ar[r]^-{d_0} & \cdots & \cdots\ar[r]^-{d_0}& X_0 & \ar[l]_-{d_1}  X_1
}
\]
 is a weak homotopy equivalence.
\end{itemize}
Examples: the nerve $N\sC$ of a small category $\sC$, which is a simplicial set and as such a simplicial space,
is a Segal space. (Note in passing that the set of $r$-simplices of $N\sC$ is the set of \emph{contravariant} functors from $I_r$ to $\sC$,
where $I_r$ is $\{0,1,\dots,r\}$ viewed as a poset with the usual ordering.) The nerve of a small topological category
$\sC$ is a Segal space provided that at least one of the maps \emph{source, target} from the space of morphisms of $\sC$ to
the space of objects of $\sC$ is a (Kan) fibration. \newline
For a Segal space $X$ it is sensible to view $X_0$ as the space of objects of something slightly more general
than a topological category, and to view $X_1$ as the space of morphisms, and $d_0,d_1\co X_1\to X_0$ as the
operators \emph{source} and \emph{target}, respectively.
In the same spirit, let
\[ \mor^h_X(a,b):=\hofiber_{(a,b)}\left[\,(d_0,d_1)\co X_1\to X_0\times X_0\,\right]  \]
for $a,b\in X_0$. There is a composition map $\mor^h_X(b,c)\times \mor^h_X(a,b) \to \mor^h_X(a,c)$, well defined at
least up to homotopy. \newline
A simplicial map $f\co X\to Y$ between simplicial spaces which satisfy $(\sigma)$ is a \emph{Dwyer-Kan} equivalence if
\begin{itemize}
\item[-] for every $a,b\in X_0$ the map
\[  \mor^h_X(a,b) \lra \mor^h_Y(f(a),f(b)) \]
induced by $f$ is a weak equivalence;
\item[-] \emph{essential surjectivity}: for every $c\in Y_0$ there exist $b\in X_0$ and an element of $\mor^h_Y(f(b),c)$ which is
weakly invertible.
\end{itemize}
A simplicial space $X$ is a \emph{complete Segal space} if in addition to $(\sigma)$ it satisfies
another property $(\kappa)$, completeness, for which the reader can consult \cite[\S2]{BoavidaWeissLong}.
The point of this additional condition is the following. Firstly, every Segal space admits a Dwyer-Kan
equivalence to a complete Segal space. Secondly, a Dwyer-Kan equivalence between two Segal spaces which are both complete
is a degreewise weak equivalence.

Rezk's model category structure on the category $s\sS$ of simplicial spaces can be described roughly as follows.
Start with the standard model category structure on $\sS$, the category of spaces (=simplicial sets). Use it to define a
model category structure on $s\sS$ where either the weak equivalences and the cofibrations
are defined levelwise, or the weak equivalences and the fibrations are defined levelwise (yes, there are two options).
Write $s\sS_1$ for $s\sS$ with this model category structure. There is a unique model category structure on $s\sS$
which has the same cofibrations as $s\sS_1$ and which for the fibrant objects has the complete Segal spaces
which are also fibrant as objects of $s\sS_1$. See \cite[\S B.2]{BoavidaWeissLong}. This is the Rezk model category structure.
Write $s\sS_2$ for $s\sS$ with that model category structure. A map $f\co X\to Y$ of simplicial spaces is a
weak equivalence in $s\sS_2$ if and only if, for every complete Segal space $Z$ which is fibrant in $s\sS_1$\,, the induced map
\[  \rmap_{s\sS_1}(Y,Z) \to \rmap_{s\sS_1}(X,Z) \]
(of derived mapping spaces formed in $s\sS_1$) is a weak equivalence. This characterization of the
weak equivalences means that $s\sS_2$ can also be constructed from $s\sS_1$ by a left localization.

We will also need a model category structure on the over category $s\sS/Z$ where $Z$ is an object of $s\sS$, for us typically
a Segal space but not a complete one. There is a standard way by which a model category structure on some category $\sC$
determines a model category structure on each of the over categories $\sC/c$\,, for objects $c$ of $\sC$.
See \cite[Expl.1.7]{GoerssSchemmerhorn}. We apply this with $s\sS_2$ and we denote the result by $s\sS_2/Z$. By definition the fibrant objects in $s\sS_2/Z$ are simply the fibrations with target $Z$ in $s\sS_2$. But they also
have a characterization as the morphisms $Y\to Z$ in $s\sS$ which are fibrations in $s\sS_1$ and make $Y$ into a
\emph{fiberwise complete Segal space over $Z$};
see \cite[\S B]{BoavidaWeissLong} for more details.

\medskip
Cisinski and Moerdijk \cite{CisinskiMoerdijk2} have an analogous framework for topological operads; see also related
earlier papers by Cisinski and Moerdijk \cite{CisinskiMoerdijk1}, 
and by Moerdijk and I.Weiss \cite{MoerdijkWeiss}. The starting point for this is a small category $\dend$ (in their notation, $\Omega$)
whose objects are certain finite trees. In more detail, an object $T$ of $\dend$ is a finite nonempty set $\epsilon(T)$ with a
partial order $\le$ and a distinguished subset $\lambda(T)$ of the set of maximal elements of $\epsilon(T)$ such that the following conditions are
satisfied:
\begin{itemize}
\item $\epsilon(T)$ has a minimal element (called the \emph{root});
\item for each element $e$ of $\epsilon(T)$, the set $\{y\in \epsilon(T)~|~y\le e\}$ with the restricted ordering is linearly ordered.
\end{itemize}
The elements of $\epsilon(T)$ are also called \emph{edges} of the tree $T$. The elements of $\lambda(T)$ are called the \emph{leaves}
of $T$. (The elements of $\epsilon(T)\smin\lambda(T)$
are sometimes called \emph{vertices}, but I have found this confusing and I prefer to call them \emph{non-leaf edges}.)
An object $T$ of $\dend$ generates a finite (colored!) operad with color set $\epsilon(T)$;
for each non-leaf edge $a$ in $T$ there is a generating operation with target $a$ and with multi-source equal to the set of
edges which are just above $a$ in the ordering. A morphism $S\to T$ in $\dend$ is by definition a morphism of the associated
finite operads. See \cite[\S7.2]{BoavidaWeissLong} for more details and examples. The category $\dend$ contains a copy of the
category $\Delta$ (with objects $[k]=\{0,1,\dots,k\}$ for $k\ge 0$ and with monotone maps as morphisms). A \emph{dendroidal space}
is a functor from $\dend^\op$ to spaces. Therefore a dendroidal space determines a simplicial space by restriction from $\dend^\op$ to
$\Delta^\op$. 
\newline
There is a concept of \emph{Segal dendroidal space} analogous to the concept of Segal space, and
a corresponding model category structure on the category $d\sS$ of dendroidal spaces. Again, this can be obtained by starting with a standard model
category structure on $d\sS$, for which we write $d\sS_1$~, and then applying a localization process. The result is $d\sS_2$ which has the
same underlying category and the same cofibrations as $d\sS_1$~, but the fibrant objects in $d\sS_2$ are the Segal dendroidal spaces
which are also fibrant in $d\sS_1$. \newline
A topological operad $Q$ has a dendroidal nerve $N_dQ$, which is a dendroidal space such that $(N_dQ)_T$ is the
space of operad maps from the operad associated with the tree $T$ to $Q$. If $Q$ is monochromatic, then $(N_dQ)_T$ is a point 
for $T=\eta$, the tree with a single edge. (The notation $N_dQ$ for the dendroidal
nerve of an operad $Q$ is widely used. It clashes with the notation $N_d\sC$ for the set of $d$-simplices in the nerve of a 
category $\sC$. Notation like $(N\sC)_d$ can be used in such cases.)
The dendroidal nerve $N_dQ$ is then a Segal dendroidal space. (It does not make sense
to insist on a completeness property as in \cite{Rezk}, or to claim such a property, if we want to work with monochromatic operads. I 
am indebted to the referee for drawing my attention to this point. By analogy, the bar construction alias nerve of a topological monoid $M$ is a 
simplicial space which qualifies as a Segal space, but it is a \emph{complete} Segal space only if the subspace of weakly invertible elements in $M$ 
is weakly contractible. See also \cite[Expl.7.7]{BoavidaWeissLong}. Similarly, 
if the monochromatic operad $Q$ satisfies $Q(1)\simeq *$, then $N_dQ$ can call itself a \emph{complete} 
Segal dendroidal space, but we are not imposing that condition yet.) A key result of the dendroidal theory is that the dendroidal nerve functor induces a weak equivalence
\begin{equation} \label{eqn-dnerveconversion}
\xymatrix{
\rmap(P,Q) \ar[r]^-\simeq & \rmap(N_dP,N_dQ)
}
\end{equation}
when $P$ and $Q$ are (monochromatic) operads. 
In the left-hand side, the derived mapping space $\rmap(P,Q)$ can be interpreted using a standard
model structure with levelwise weak equivalences and fibrations. In the right-hand side, use either
$d\sS_1$ or $d\sS_2$~; it does not matter which because $N_dQ$ is a Segal dendroidal space.

\medskip
For a simplicial set $Y$ let $\simp(Y)$ be the category whose objects are pairs $(m,y)$ with $y\in X_m$~; a morphism from
$(m,y)$ to $(n,z)$ is a morphism $f\co [m]\to [n]$ in $\Delta$ such that $f^*(z)=y$. There is a functor
\[  \varphi\co\simp(N\fin) \lra \dend \]
defined as follows. To an object $(p,S_*)$ of $\simp(N\fin)$, where
\[  S_*=(S_0\leftarrow S_1\leftarrow S_2 \leftarrow \cdots\leftarrow S_p) \]
associate the tree $T$ where $\epsilon(T)$ is the disjoint union of the $S_i$ and an additional element $r$, with $\lambda(T)$
corresponding to $S_p$\,.
The partial order on $\epsilon(T)$ is the obvious one: $r$ is the minimal element and $x\in S_i$ is
$\le y\in S_j$ if $i\le j$ and the composite map from $S_j$ to $S_i$ in the string $S_*$ takes $y$ to $x$.
The functor $\varphi$ establishes up a close relationship between dendroidal spaces and simplicial spaces over $N\fin$. Indeed,
a simplicial space over $N\fin$ is the same thing as a contravariant functor from $\simp(N\fin)$ to spaces. Therefore
pre-composition with $\varphi$ is a functor $\varphi^*$ which takes us from dendroidal spaces to simplicial spaces over
$N\fin$. \emph{Important example}: for a topological operad there is an isomorphism of $\varphi^*(N_dP)$
with $N\sC_P$~, both viewed as simplicial spaces over $N\fin$.

\medskip
The main result of \cite[\S7]{BoavidaWeissLong} is that for topological operads $P$ and $Q$ the functor
$\varphi^*$ induces a weak equivalence
\begin{equation} \label{eqn-dendfin}
\xymatrix@R=10pt{
\rmap(N_dP,N_dQ) \ar[r] & \rmap_{N\fin}(\varphi^*N_dP,\varphi^*N_dQ) \ar@{=}[d]  \\
& \rmap_{N\fin}(N\sC_P,N\sC_Q)
}
\end{equation}
provided $P(\uli 0),P(\uli 1),Q(\uli 0)$ and $Q(\uli 1)$ are weakly contractible. In the right-hand side
we can use the model category structure $s\sS_1/N\fin$ or $s\sS_2/N\fin$~; it does not matter because
$N\sC_Q$ is a fiberwise complete Segal space over $N\fin$. --- In view of
(\ref{eqn-dnerveconversion}) this implies that the nerve functor determines a weak equivalence
\begin{equation} \label{eqn-Nfinconversion}
\xymatrix{
\rmap(P,Q) \ar[r]^-\simeq & \rmap_{N\fin}(N\sC_P,N\sC_Q)
}
\end{equation}
under the same conditions on $P$ and $Q$.

\begin{expl} {\rm Suppose that $P$ is the operad of little $m$-disks and $Q$ is the operad of little $n$-disks.
Then $N\sC_P$ and $N\sC_Q$ are weakly equivalent, over $N\fin$, to certain simplicial spaces
$\config(\RR^m)$ and $\config(\RR^n)$ defined as nerves of certain (topological) categories of configurations
in $\RR^m$ and $\RR^n$, respectively. See \cite[\S3]{BoavidaWeissLong} and \cite[expl.7.2]{BoavidaWeissLong}.
With Andrade's particle models for the configuration categories, $M\mapsto \config(M)$ is a continuous functor on the category of topological
manifolds and injective continuous maps. It follows that there are compatible actions of the homeomorphism group of $\RR^m$ (from the right)
and the homeomorphism group of $\RR^n$ (from the left) on
\[  \rmap_{N\fin}(\config(\RR^m),\config(\RR^n)) \simeq \rmap_{N\fin}(N\sC_P,N\sC_Q). \]
Similarly, there is an interesting map from the space of injective continuous maps $\RR^m\to \RR^n$ to the space
$\rmap_{N\fin}(\config(\RR^m),\config(\RR^n))$. All these good features are not easy to see
in the operadic description $\rmap(P,Q)$.
}
\end{expl}

\section{Truncated operads and truncated dendroidal spaces}
For an integer $k\ge 1$, a $k$-truncated operad in the symmetric monoidal category of spaces
is defined like an operad in spaces except for the following changes: all operations have arity $\le k$
and composition of operations is only defined where it does not contradict this restriction.
More specifically, we can describe a $k$-truncated (monochromatic, symmetric) operad in spaces in the following terms. It is
a functor $P$ from the category of finite sets of cardinality $\le k$ and bijections to
spaces, and for every map $f:T\to S$ of finite sets of cardinality $\le k$ there is an operation
\[  \lambda_f\co P(S)\times\prod_{i\in S} P(T_i) \lra P(T) \]
where $T_i=f^{-1}(i)\subset T$. Also $P(S)$ contains a distinguished unit element when $S$ is a singleton.
Sensible naturality, associativity and unital properties are satisfied. In particular any
permutation $f\co S\to S$, where $|S|\le k$, induces a map $P(S)\to P(S)$ in two ways: firstly because
$P$ is a functor, and secondly by
\[ P(S)\ni x\mapsto \lambda_f(x,1,1,\dots,1)\in P(S)~. \]\
for $x\in P(S)$. These two maps agree as per definition.

\medskip
For an integer $k\ge 1$ let $\dend_k\subset \dend$ be the full subcategory whose objects are the
trees $T$ such that for every $t\in \epsilon(T)\smin\lambda(T)$ the set $\{s\in \epsilon(T)~|~s>t\}$ has at most
$k$ minimal elements. (These minimal elements are often called the \emph{incoming edges} ... of a fictional vertex associated with the non-leaf edge $t$.)
Similarly let $\fin_{\le k}$ be the full subcategory spanned by the objects $\uli m$ where $m\le k$.
The functor $\varphi$ above restricts to a functor
\[  \varphi_k\co \simp(N\fin_{\le k}) \lra \dend_k\,. \]
An object $T$ of $\dend_k$ generates a finite $k$-truncated
(colored) operad with color set $\epsilon(T)$; for each non-leaf edge $a$ in $T$ there is a generating operation
with target $a$ and with multi-source equal to the set of edges which are just above $a$ in the ordering.
A $k$-truncated operad $P$ has a nerve $N_dP$ which we regard as a
contravariant functor from $\dend_k$ to spaces. It is defined in such a way that $(N_dP)_T$
is the space of $k$-truncated operad morphisms from the $k$-truncated operad
associated with $T$ to $P$. The $k$-truncated operad $P$ also determines a topological category
$\sC_P$ (category in spaces) with a forgetful functor to $\fin_{\le k}$. The space of objects of $\sC_P$ is
\[  \coprod_{m=0}^k P(\uli m). \]
For a morphism $f\co \uli\ell\to \uli m$ in $\fin_{\le k}$ the space of morphisms in $\sC_P$ lifting $f$
is
\[  P(\uli\ell) \times \prod_{j\in\uli m} P(\uli\ell_{\,j}) \]
where $\ell_j$ is the cardinality of $f^{-1}(j)$. \newline
Keeping in mind that a simplicial space over $N\fin_{\le k}$ is the same thing as a
contravariant functor from $\simp(N\fin_{\le k})$, we can say that composition with $\varphi_k$
is a functor $\varphi^*_k$ from $k$-truncated dendroidal spaces to simplicial spaces over $N\fin_{\le k}$\,.
In that sense there is an isomorphism
\[  \varphi_k^*(N_dP)\cong N\sC_P \]
of simplicial spaces over $N\fin_{\le k}$\,, assuming that $P$ is a $k$-truncated operad in spaces.

\begin{thm} \label{thm-opersuppl} Let $P$ and $Q$ be $k$-truncated operads (in spaces, $k\ge 1$) for which the spaces
$P(\uli 0)$, $P(\uli 1)$, $Q(\uli 0)$ and $Q(\uli 1)$ are weakly contractible. Then composition with the
functor $\varphi_k^*$ is a weak equivalence
\[
\xymatrix@R=10pt{
\rmap(N_dP,N_dQ) \ar[r]^-\simeq & \rmap_{N\fin_{\le k}}(\varphi_k^*N_dP,\varphi^*N_dQ) \ar@{=}[d]  \\
& \rmap_{N\fin_{\le k}}(N\sC_P,N\sC_Q).
}
\]
\end{thm}

This will be proved in the following sections. For clarification, we make sense of $\rmap(N_dP,N_dQ)$ using a
model category structure with levelwise weak equivalences on the category of contravariant
functors from $\dend_k$ to spaces. We make sense of $\rmap_{N\fin_{\le k}}(N\sC_P,N\sC_Q)$ using a model category structure
on the category of simplicial spaces over $N\fin_{\le k}$ with levelwise weak equivalences.
(See \cite{DwyKa2}, where it is shown that the derived mapping spaces in a model category depend
mainly on the subcategory of weak equivalences, not much on the subcategories of cofibrations and fibrations,
respectively.) The proof of theorem~\ref{thm-opersuppl} given here works equally well in the untruncated setting, $k=\infty$.
It can be seen as another way to show that~(\ref{eqn-dendfin}) is a weak equivalence which happens to generalize
easily to the truncated situation.

\section{Leaves are unnecessary}
In \cite[\S7]{BoavidaWeissLong} the weak equivalence~(\ref{eqn-dendfin}) is established in essentially two steps.
In the first step, which is the easier one, the category $\dend$ gets replaced by a much more accessible
subcategory $\dendroca$. It is the subcategory of $\dend$ obtained by allowing as objects only the trees $T$ with
empty leaf set $\lambda(T)$ and as morphisms only the morphisms in $\dend$ between trees with no leaves which take root to root. \newline
An object $T$ of $\dendroca$ can therefore be described as a finite set $T$ with an order relation (written $\le$) such that
\begin{itemize}
\item[-] $T$ has a (unique) minimal element, called the \emph{root};
\item[-] for every $s\in T$ the set $\{t\in T~|~t\le s\}$ is linearly ordered (with the ordering
induced from $T$).
\end{itemize}
(We need not distinguish anymore between $T$ and the edge set of $T$.)
A morphism $S\to T$ in $\dendroca$ is a map $f$ from $S$ to $T$ which respects the order relations, takes root
to root and has the additional property that it preserves \emph{independence}. That is, if for $u,v\in S$ we have neither
$u\leq v$ nor $v\leq u$, then for $f(u),f(v)\in T$ we have neither $f(u)\leq f(v)$
nor $f(v)\leq f(u)$.

There is a functor $\psi\co \simp(N\fin)\to \dendroca$ very similar to $\varphi\co \simp(N\fin)\to \dend$.
For an object $(p,S_*)$ of $\simp(N\fin)$, where
\[  S_*=(S_0\leftarrow S_1\leftarrow S_2 \leftarrow \cdots\leftarrow S_p), \]
let $\psi(p,S_*)$ be tree $T$ without leaves which, as a set, is the disjoint union of the $S_i$ and an additional element $r$.
The partial order on $T$ is the obvious one: $r$ is the minimal element and $x\in S_i$ is
$\le y\in S_j$ if $i\le j$ and the composite map from $S_j$ to $S_i$ in the string $S_*$ takes $y$ to $x$.

The formal relationship between $\varphi$ and $\psi$ is a little more complicated than one might expect.
The inclusion $\iota\co \dendroca\to \dend$ has a left adjoint $\kappa\co \dend\to \dendroca$. Equations such
as $\varphi=\iota\psi$ or $\psi=\kappa\varphi$ come to mind, but both are false. Instead we have
$\varphi\beta=\iota\psi$~,
where $\beta\co \simp(N\fin) \to \simp(N\fin)$ is the endofunctor which takes a string
\[  (S_0\leftarrow S_1\leftarrow S_2 \leftarrow \cdots\leftarrow S_p) \]
to the string
\[ (S_0\leftarrow S_1\leftarrow S_2 \leftarrow \cdots\leftarrow S_p\leftarrow \emptyset). \]
Now let $N_d^{rc}P$ and $N_d^{rc}Q$ be the restrictions of $N_dP$ and $N_dQ$,
respectively, to $\dendroca$. On the basis of the observations just above, showing
that~(\ref{eqn-dendfin}) is a weak equivalence reduces easily to showing that the map
\begin{equation}  \label{eqn-denfinnrc}
\rmap(N_d^{rc}P,N_d^{rc}Q) \lra \rmap(\psi^*N_d^{rc}P,\psi^*N_d^{rc}Q) \end{equation}
obtained by composition with $\psi^*$ is a weak equivalence. We are still assuming that $P(\uli 0),P(\uli 1)$ and
$Q(\uli 0)$, $Q(\uli 1)$ are weakly contractible. 
(There is a commutative diagram
\[
\xymatrix@C=50pt{
\rmap(\varphi^*N_dP,\varphi^*N_dQ) \ar[d]^-{\beta^*} & \ar[l]_-{(\ref{eqn-dendfin})} \rmap(N_dP,N_dQ) \ar[d]^-{\iota^*} \\
\rmap(\psi^*N_d^{rc}P,\psi^*N_d^{rc}Q) & \ar[l]_-{(\ref{eqn-denfinnrc})} \rmap(N_d^{rc}P,N_d^{rc}Q)
}
\]
where the vertical arrows are weak equivalences; use $\varphi\beta=\iota\psi$ and $N_d^{rc}=\iota^*N_d$\,.)

\medskip
The message of this short section is that the same argument applies in the truncated situation.
There is a functor
\[ \psi_k\co \simp(N\fin_{\le k})\lra \dend_k\cap \dendroca \]
obtained by restriction of $\psi$. Suppose that $P$ and $Q$ are $k$-truncated operads
for which $P(\uli 0),P(\uli 1)$ and
$Q(\uli 0)$, $Q(\uli 1)$ are weakly contractible. Let $N_d^{rc}P$ and $N_d^{rc}Q$
be the contravariant functors from $\dend_k\cap \dendroca$ to spaces obtained by
restricting $N_dP$ and $N_dQ$, respectively. Then for
the proof of theorem~\ref{thm-opersuppl} it suffices to show that the map
\begin{equation}  \label{eqn-denfinnrctrunc}
\rmap(N_d^{rc}P,N_d^{rc}Q) \lra \rmap(\psi_k^*N_d^{rc}P,\psi_k^*N_d^{rc}Q) \end{equation}
obtained by composing with $\psi_k^*$ is a weak equivalence.

\section{Bridging the gap}
Let $\RR\psi_*$ be the homotopy right Kan extension along $\psi$. This is applicable to contravariant functors
from $\simp(N\fin)$ to spaces and yields contravariant functors from $\dendroca$ to spaces.
It serves as a homotopy right adjoint to the functor $\psi^*$ given by pre-composition with $\psi$.
It is shown in \cite[\S7]{BoavidaWeissLong} that, under conditions on $Q$ as in~(\ref{eqn-denfinnrc}),
the homotopy unit
\begin{equation} \label{eqn-hounit} N_d^{rc}Q \lra \RR\psi_*\psi^*(N_d^{rc}Q) \end{equation}
is a weak equivalence. This implies in a formal manner that the map~(\ref{eqn-denfinnrc})
is a weak equivalence. See \cite[lem.A.1]{BoavidaWeissLong}.

This type of argument is also available in the truncated setting, but showing that
the truncated analogue of~(\ref{eqn-hounit}) is a weak equivalence is harder than
showing that~(\ref{eqn-hounit}) is a weak equivalence. Therefore we proceed in two steps by writing
the functor $\psi$ and its truncated variant $\psi_k$ as a composition of two functors. In doing so
we deviate from the line of reasoning developed in \cite[\S7]{BoavidaWeissLong}. It amounts to additional work,
but there is the surprising reward that we can avoid the use of a difficult lemma \cite[lem.7.14]{BoavidaWeissLong}.

\begin{defn} {\rm There is a category $\dendlev$ which is halfway between $\simp(N\fin)$ and
$\dendroca$. An object of $\dendlev$ is an object of $\simp(N\fin)$. A morphism in $\dendlev$
from $S_*=(S_0\leftarrow S_1\leftarrow\cdots\leftarrow S_k)$
to $R_*=(R_0\leftarrow R_1\leftarrow\cdots\leftarrow R_\ell)$ consists of a monotone map
$u\co [k]\to[\ell]$ and monotone injections $v_j\co S_j\to R_{u(j)}$, one for every $j\in [k]$, such that
the diagrams
\[
\xymatrix{S_{j-1} \ar[d]^{v_{j-1}} && \ar[ll] \ar[d]^{v_j} S_j \\
R_{u(j-1)} & \ar[l] \cdots  & \ar[l] R_{u(j)}
}
\]
commute for $j\in \{1,\dots,k\}$. If the $v_j$ are bijective, then they are necessarily
identity maps and the collection $(u,(v_j))$ is a morphism in $\simp(N\fin)$ from $S_*$ to $R_*$\,.
Therefore $\simp(N\fin)\subset \dendlev$.
}
\end{defn}

Let us note that in the category $\simp(N\fin)$ or, for that matter, in any category of the form $\simp(Y)$ where $Y$ is a
simplicial set, no object admits nontrivial automorphisms. The same can be said of $\dendlev$: no object
admits nontrivial automorphisms.

\smallskip
Write $\alpha\co \simp(N\fin)\to \dendlev$ for the inclusion.
The functor $\psi$
has an obvious extension to a functor $\xi\co \dendlev \to \dendroca$, so that there is a commutative
triangle
\begin{equation} \label{eqn-decomp}
\begin{split}
\xymatrix@C=7pt{ \ar[dr]_-{\alpha} \simp(N\fin) \ar[rr]^-{\psi} && \dendroca \\
& \dendlev  \ar[ur]_-{\xi}
}
\end{split}
\end{equation}

\begin{lem} \label{lem-alpha} Let $Z=\xi^*(N_d^{rc}Q)$. The homotopy unit
$Z \to \RR\alpha_*\alpha^*Z$
is a weak equivalence.
\end{lem}

\proof
Let $R_*=(R_0\leftarrow R_1\leftarrow\cdots\leftarrow R_\ell)$ be an object of $\dendlev$\,.
For $(\RR\alpha_*\alpha^*Z)(R_*)$ we have the standard formula
\[  \holimsub{S_*\to R_*} Z(S_*) \]
where the homotopy limit is taken over a certain comma category or over category $\mathscr U$ which depends on $R_*$.
The objects of $\mathscr U$ are pairs consisting of an object $S_*$ of $\simp(N\fin)$ and a morphism
$f\co S_*\to R_*$ in $\dendlev$. A morphism in $\mathscr U$ from $(S_*,f)$ to $(S'_*,g)$
is a morphism $S_*\to S'_*$ in $\simp(N\fin)$ which is over $R_*$ when viewed as a morphism in $\dendlev$. \newline
We introduce full subcategories
\[  \mathscr U_{-1}\supset \mathscr U_0\supset \mathscr U_1\supset \mathscr U_2 \supset\cdots \supset
\mathscr U_{\ell-1} \supset \mathscr U_\ell \]
of $\mathscr U$, where $\mathscr U_{-1}=\mathscr U$ and $\mathscr U_\ell$ is the comma category determined by the identity functor
$\simp(N\fin)\to\simp(N\fin)$ and the object $R_*$\,. (The integer $\ell$ is determined by $R_*$\,.)
The details are as follows. An object of $\mathscr U$ given by $S_*\to R_*$~, or more precisely, by a commutative diagram
\begin{equation} \label{eqn-Uladder}
\begin{split}
\xymatrix@C=15pt{\cdots && \ar[ll] S_{j-1} \ar[d]^{v_{j-1}} && \ar[ll] \ar[d]^{v_j} S_j && \ar[ll] \cdots \\
& \ar[l] \cdots & \ar[l] R_{u(j-1)} & \ar[l] \cdots  & \ar[l] R_{u(j)} & \ar[l] \cdots & \ar[l]
}
\end{split}
\end{equation}
belongs to $\mathscr U_r$ if and only if $v_j$ is bijective (hence an identity map) for all $j$ such that $u(j)\le r$.
In particular that object $S_*\to R_*$ belongs to $\mathscr U_\ell$ if and only if $S_*\to R_*$ is a morphism in $\simp(N\fin)$.
We use the abbreviation $\bar Z(S_*\to R_*):= Z(S_*)$ for an object $S_*\to R_*$ in $\mathscr U$.
Then $\bar Z$ is a contravariant functor from $\mathscr U$ to spaces and
\[ (\RR\alpha_*\alpha^*Z)(R_*)=\holim~\bar Z=\holim~\bar Z|\mathscr U_{-1}\,. \]
There is a string of forgetful projections
\[ \holim~\bar Z|\mathscr U_{-1} \lra \holim~\bar Z|\mathscr U_{0} \lra \cdots \lra \holim~\bar Z|\mathscr U_{\ell}\,. \]
We have our unit map from
$Z(R_*)$ to $(\RR\alpha_*\alpha^*Z)(R_*)=\holim~\bar Z|\mathscr U_{-1}$ such that the composition
\[ Z(R_*) \to \holim~\bar Z|\mathscr U_{-1} \lra \holim~\bar Z|\mathscr U_{0} \lra \cdots \lra \holim~\bar Z|\mathscr U_{\ell} \]
is a weak equivalence for rather trivial reasons. Therefore it suffices to show that
each of the projection maps $\holim~\bar Z|\mathscr U_r\to \holim~\bar Z|\mathscr U_{r+1}$ admits a homotopy left inverse,
making $\holim~\bar Z|\mathscr U_r$ a
homotopy retract of $\holim~\bar Z|\mathscr U_{r+1}$. To achieve that we shall construct two functors
\[  V\co \mathscr U_r\to \mathscr U_r~,\qquad W\co \mathscr U_r \to \mathscr U_{r+1} \]
and natural transformations $\id\Rightarrow V\Leftarrow W$,
where the functor $V$ takes $\mathscr U_{r+1}$ to itself. The crucial property is that, for every
object $S_*\to R_*$ in $\mathscr U_r$~, the natural morphism $W(S_*\to R_*)\to V(S_*\to R_*)$ is taken to a weak
equivalence by the functor $\bar Z$. Then we shall have the following maps:
\[
\xymatrix{
\holim~\bar Z|\mathscr U_{r+1} \ar[r] & \holim~(\bar Z|\mathscr U_{r+1})\circ W \\
& \ar[u]^\simeq \holim~(\bar Z|\mathscr U_r)\circ V \ar[r] & \holim~\bar Z|\mathscr U_r
}
\]
(the first by prolongation, the other two using the natural transformations) which give us the required homotopy class of maps
from $\holim~\bar Z|\mathscr U_{r+1}$ to $\holim~\bar Z|\mathscr U_{r}$\,. \newline
For the description of $V$, imagine an object $S_*\to R_*$ of $\mathscr U_r$ given by a diagram like~(\ref{eqn-Uladder}),
where $j$ runs through $\{0,1,\dots,k\}$.
Determine the unique $t\in \{0,1,\dots,k,k+1\}$ such that $u(j)\le r$ whenever $j<t$ and $u(j)>r$ whenever $j\ge t$.
Let $V(S_*\to R_*)=(S'_*\to R_*)$ where
\[  S'_*=(S_0\leftarrow S_1\leftarrow \cdots \leftarrow S_{t-1}\leftarrow R_{r+1}
\leftarrow S_t\leftarrow S_{t+1} \leftarrow\cdots \leftarrow S_k), \]
that is, $S'_j=S_j$ for $j<t$, $S'_t=R_{r+1}$ and $S'_j=S_{j-1}$ for $t<j\le k+1$. The arrow
from $S'_{t+1}=S_t$ to $S'_t=R_{r+1}$ is $v_t$\,. The morphism (in $\dendlev$) from $S'_*$ to $R_*$ is defined in such a way
that there is a commutative triangle
\[
\xymatrix@C=8pt@R=12pt{
S_* \ar[rr] \ar[dr] && \ar[dl] S'_* \\
& R_*
}
\]
in $\dendlev$\,, where the horizontal arrow (a morphism in $\simp(N\fin)$) is defined by the monotone
injection $[k]\to[k+1]$ which omits $t$. Of course, $S'_t$ should be mapped to $R_{r+1}=R_{u(t)}$ by the identity.
This commutative triangle contributes to a sketchy description not only of $V$, but
also of our preferred natural transformation from $\id$ to $V$. Now the functor $W$ is defined, on an object $S_*\to R_*$ of
$\mathscr U$ as before, by starting from $V(S_*\to R_*)=(S'_*\to R_*)$ as described and erasing from $S'_*$ all the terms
$S'_{j+1}=S_j$ where $j\ge t$ and $v(j)=r+1$. Call the result $S''_*\to R_*$\,. Again there is a
commutative triangle
\[
\xymatrix@C=8pt@R=12pt{
S'_* \ar[dr] && \ar[ll] \ar[dl] S''_* \\
& R_*
}
\]
where the horizontal arrow (a morphism in $\simp(N\fin)$) is defined by a monotone injection $[k+1-a]\to[k+1]$ which omits an interval of
the form $\{t+1,\dots,t+a\}$. This contributes to a sketchy description not only of $W$, but also
of our preferred natural transformation from $W$ to $V$. The morphism in $\mathscr U_r$
defined by this triangle is taken to a homotopy equivalence by the functor $\bar Z$. To illustrate that, here is a picture
describing $S'_*$ and $S''_*$ and the preferred morphism
between them:
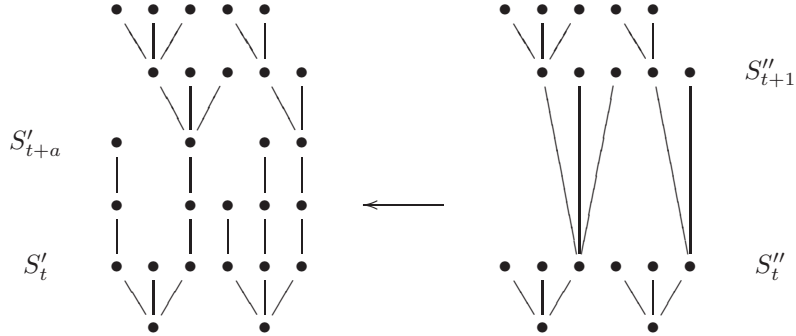
\[
\xymatrix@C=3pt@R=2pt{
&&\bullet\ar@{-}[ddr] & \bullet\ar@{-}[dd] & \bullet\ar@{-}[ddl] & \bullet\ar@{-}[ddr] & \bullet\ar@{-}[dd]  &&&&&&&&&
\bullet\ar@{-}[ddr] & \bullet\ar@{-}[dd] & \bullet\ar@{-}[ddl] & \bullet\ar@{-}[ddr] & \bullet\ar@{-}[dd]
\\
\\
&&& \bullet\ar@{-}[ddr] & \bullet\ar@{-}[dd] & \bullet\ar@{-}[ddl] & \bullet\ar@{-}[ddr] & \bullet\ar@{-}[dd] &&&&&&&&
& \bullet\ar@{-}[ddddddr] & \bullet\ar@{-}[dddddd] & \bullet\ar@{-}[ddddddl] & \bullet\ar@{-}[ddddddr] & \bullet\ar@{-}[dddddd] && S''_{t+1} \\
\\
S'_{t+a} && \bullet\ar@{-}[dd] &&  \bullet \ar@{-}[dd] &   & \bullet\ar@{-}[dd]
& \ar@{-}[dd] \bullet &&&&&&&&&&&&& \\
 &&&&&&&&&&&&&&&&&&&&&&  \\
&& \bullet\ar@{-}[dd] &&  \bullet \ar@{-}[dd] & \bullet  \ar@{-}[dd]  & \bullet \ar@{-}[dd]
& \ar@{-}[dd] \bullet &&&&&&\ar[llll]&&&&&&& \\
\\
S'_t && \bullet\ar@{-}[ddr] & \bullet \ar@{-}[dd] &  \bullet \ar@{-}[ddl] & \bullet  \ar@{-}[ddr]  & \bullet \ar@{-}[dd]
& \ar@{-}[ddl] \bullet &&&&&&&& \bullet\ar@{-}[ddr] & \bullet \ar@{-}[dd] &  \bullet \ar@{-}[ddl] & \bullet  \ar@{-}[ddr]  & \bullet \ar@{-}[dd]
& \ar@{-}[ddl] \bullet && S''_t \\
&&&&&&&&&&&&&&&&&&&&&&  \\
&&& \bullet  &&& \bullet  &&&&&&&&&& \bullet  &&& \bullet  &
}
\]
We are saying that the corresponding operator in $N_d^{rc}Q$ is a weak equivalence.
This is based on the assumption that $Q(\uli 0)$ and $Q(\uli 1)$ are weakly contractible. \qed

\begin{lem} \label{lem-xi} Let $Z=N_d^{rc}Q$. The homotopy unit
$Z \to \RR\xi_*\xi^*Z$
is a weak equivalence.
\end{lem}

\proof
Let $T$ be an object of $\dendroca$. The standard
formula for $(\RR\xi_*\xi^*Z)(T)$ is
\[   \holimsub{\xi(S_*)\to T} Z(\xi(S_*))~. \]
Here the homotopy inverse limit is taken over a comma category $\mathscr V(T)$. An object of $\mathscr V(T)$
is an object $S_*$ of $\dendlev$ together with a morphism $f\co \xi(S_*)\to T$ in $\dendroca$\,. A morphism in $\mathscr V(T)$
from $(S_*,f)$ to $(R_*,g)$ is a morphism $S_*\to R_*$ in $\dendlev$ which turns into a morphism
over $T$ on applying $\xi$. Let $F_T$ from $\mathscr V(T)$ to spaces be the functor which takes
$(S_*,f)$ to $Z(\xi(S_*))$. Then we can write
\[  (\RR\xi_*\xi^*Z)(T)= \holim~F_T~. \]
We proceed by induction on the number of nodes of $T$, where \emph{node} means a vertex with more than one incoming edge.
The induction beginning includes the case where $T$ has zero nodes. Then $T$ is linearly ordered. It is easy to see
that $Z(T)$ is weakly contractible, since we are assuming weak contractibility of $Q(\uli 0)$ and $Q(\uli 1)$.
Also, for each object $(S_*,f)$ in $\mathscr V(T)$, the space $Z(\xi(S_*))=F_T(S_*,f)$
is weakly contractible since $\xi(S_*)$ is linearly ordered. Therefore, if $T$ has zero nodes, the unit map
from $Z(T)$ to $(\RR\xi_*\xi^*Z)(T)$ is a weak equivalence. \newline
The induction beginning also includes the case where $T$ has exactly one node. This case is surprisingly hard. Let
$\mathscr V_0(T)$ be the full subcategory of $\mathscr V(T)$ obtained by deleting all objects $(S_*,f)$
of $\mathscr V(T)$ where $\xi(S_*)$ has zero nodes, or equivalently, the sets $S_i$ all have cardinality $\le 1$.
It is easy to see that the restriction map
\[  \holim~F_T \lra \holim~(F_T|\mathscr V_0(T)) \]
is a weak equivalence, since the value of $F_T$ on any object of $\mathscr V(T)$ not in $\mathscr V_0(T)$ is weakly contractible
and since any morphism in $\mathscr V(T)$ with source in $\mathscr V_0(T)$ has target in $\mathscr V_0(T)$.
Next, let $\mathscr V_1(T)$ be the full subcategory of $\mathscr V_0(T)$ consisting of those objects $(S_*,f)$
where the set $S_0$ has cardinality $>1$. (We write
\[  S_*=(S_0\leftarrow S_1\leftarrow \cdots\leftarrow S_{k-1}\leftarrow S_k) \]
as usual.) The inclusion functor $\mathscr V_1(T)\to \mathscr V_0(T)$ has a right adjoint.
(If $(S_*,f)$ in $\mathscr V_0(T)$ has $|S_0|=|S_1|=\cdots=|S_j|=1$ and $|S_{j+1}|>1$,
then the value of that right adjoint on $(S_*,f)$ is obtained by deleting the terms $S_0,S_1,S_2,\dots,S_j$\,.)
Moreover the counit morphisms of the adjunction are taken to a weak equivalence
by $F_T$. It follows that the restriction map
\[  \holim~(F_T|\mathscr V_0(T)) \lra \holim~(F_T|\mathscr V_1(T)) \]
is a weak equivalence. Next, let $U\subset T$ be the set of incoming edges to the unique node
of $T$. It has at least two elements by our assumption on $T$. Choose a total ordering on $U$.
Let $\mathscr V_2(T)\subset \mathscr V_1(T)$ be the subcategory of $\mathscr V_1(T)$
defined as follows. An object $(S_*,f)$ of $\mathscr V_1(T)$ qualifies as an object of $\mathscr V_2(T)$
if $f$ takes the ordered set $S_0$ to $U$ by an order-preserving bijection.
A morphism $(S_*,f)\to (R_*,g)$ in $\mathscr V_1(T)$ between such objects, given by
a monotone map $u\co [k]\to[\ell]$ and monotone injections $v_j\co S_j\to R_{u(j)}$ for $j\in [k]$,
qualifies as a morphism in $\mathscr V_2(T)$ precisely if $u(0)=0$. (In that case $v_0\co S_0\to R_0$
must be an order preserving bijection, hence an identity map.) The inclusion of
$\mathscr V_2(T)$ in $\mathscr V_1(T)$ has a left adjoint. Therefore the restriction map
\[  \holim~(F_T|\mathscr V_1(T)) \lra \holim~(F_T|\mathscr V_2(T)) \]
is a weak equivalence. Now $F_T$ takes every morphism in $\mathscr V_2(T)$ to a weak
equivalence of spaces. Moreover $\mathscr V_2(T)$ clearly has an initial object $\omega$.
Together these properties imply that the projection from $\holim~(F_T|\mathscr V_2(T))$ to $F_T(\omega)$
is a weak equivalence. Putting all that together, it follows that the projection from $\holim~F_T$ itself
to $F_T(\omega)$ is a weak equivalence. Then it follows easily by inspection that the unit map $Z(T)\to (\RR\xi_*\xi^*Z)(T)=\holim~F_T\simeq F_T(\omega)$ is a weak equivalence. This completes our discussion of the case where $T$ has exactly one node.
\newline
We come to the induction step, which is rather formal and uses a mildly sheaf theoretic argument.
Suppose that $T$ is an object of $\dendroca$ such that $T=T_0\cup T_1$ where
$T_0$ and $T_1$ are subtrees (by which, at this point, we simply mean sub-posets) of $T$. More precisely we require
that if $e\in T_0$ and $e'\in T$ with $e'\le e$, we have $e'\in T_0$\,, and similarly for $T_1$. We also require that
if $e,e'\in T$ and both are incoming edges for the same vertex, then either $e,e'$ are both in $T_0$
or both not in $T_0$~; and similarly for $T_1$.
The following picture is an example.
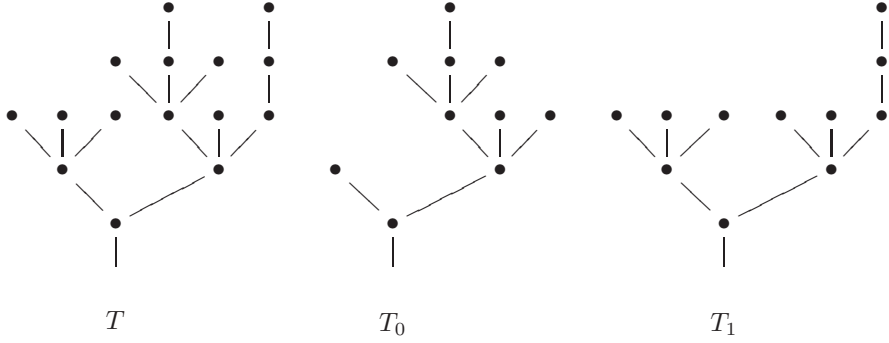
\[
\xymatrix@C=8pt@R=10pt{
&&  & \bullet \ar@{-}[d] &  & \bullet \ar@{-}[d] \\
&& \bullet \ar@{-}[dr] & \bullet \ar@{-}[d] & \bullet \ar@{-}[dl] & \bullet \ar@{-}[d] \\
\bullet \ar@{-}[dr] & \bullet \ar@{-}[d] &  \bullet \ar@{-}[dl] & \bullet \ar@{-}[dr]  & \bullet \ar@{-}[d]
& \ar@{-}[dl] \bullet \\
& \bullet \ar@{-}[dr] &&& \bullet \ar@{-}[dll] & \\
&& \bullet \ar@{-}[d] && \\
&  && \\
&& T
}
\xymatrix@C=8pt@R=10pt{
&&  & \bullet \ar@{-}[d] &  &  \\
&& \bullet \ar@{-}[dr] & \bullet \ar@{-}[d] & \bullet \ar@{-}[dl] &  \\
 &  &   & \bullet \ar@{-}[dr]   & \bullet \ar@{-}[d]
& \ar@{-}[dl] \bullet \\
& \bullet \ar@{-}[dr] &&& \bullet \ar@{-}[dll] & &\\
&& \bullet \ar@{-}[d] && \\
&  && \\
&& T_0
}
\xymatrix@C=8pt@R=10pt{
&&  &  &  & \bullet \ar@{-}[d] \\
&& &  &  & \bullet \ar@{-}[d] \\
\bullet \ar@{-}[dr] & \bullet \ar@{-}[d] &  \bullet \ar@{-}[dl] & \bullet  \ar@{-}[dr]  & \bullet \ar@{-}[d]
& \ar@{-}[dl] \bullet \\
& \bullet \ar@{-}[dr] &&& \bullet \ar@{-}[dll] & \\
&& \bullet \ar@{-}[d] && \\
&  && \\
&& T_1
}
\]
Then it
is easy to see that the square of inclusion-induced maps
\[
\xymatrix{
Z(T) \ar[r] \ar[d] & Z(T_0) \ar[d] \\
Z(T_1) \ar[r] & Z(T_0\cap T_1)
}
\]
is a homotopy pullback square. If we can show that the square of inclusion-induced maps
\[
\xymatrix{
(\RR\xi_*\xi^*Z)(T) \ar[r] \ar[d] & (\RR\xi_*\xi^*Z)(T_0) \ar[d] \\
(\RR\xi_*\xi^*Z)(T_1) \ar[r] & (\RR\xi_*\xi^*Z)(T_0\cap T_1)
}
\]
is also a homotopy pullback square, then that can pass for the induction step. With the abbreviations
or translations above, we have $(\RR\xi_*\xi^*Z)(T)=\holim~F_T$ where $F_T$ is defined on $\mathscr V_T$~.
Let $\mathscr V'_T$ be the full subcategory of $\mathscr V_T$ consisting of all objects
\[  (S_*,f\co \xi(S_*)\to T) \]
such that $\xi(S_*)$ lands in $T_0$ or in $T_1$ or even in $T_0\cap T_1$\,. By inspection, the square of restriction maps
\[
\xymatrix{
\holim~(F_T|\mathscr V'_T) \ar[r] \ar[d] & \holim~F_{T_0} \ar[d] \\
\holim~F_{T_1} \ar[r] & \holim~F_{T_0\cap T_1}
}
\]
is a homotopy pullback square. Therefore it is enough to show that the restriction map
from $\holim~F_T$ to $\holim~(F_T|\mathscr V'_T)$ is a weak equivalence. A
formula of Dwyer-Kan \cite[9.7]{DwyKa1} or Bousfield-Kan allows us to
identify $\holim~(F_T|\mathscr V'_T)$ with the homotopy inverse limit
of $\RR\eta_*(F_T|\mathscr V'_T)$~, where $\eta\co\mathscr V'_T\to \mathscr V_T$ is the inclusion
and $\RR\eta_*$ denotes the homotopy right Kan extension along $\eta$.
In the definition of $\RR\eta_*$ we use, for each object $(S_*,f)$ in $\mathscr V_T$~,
the comma category $\eta/(S_*,f)$. In that comma category there is a diagram
(of three objects)
\[
\xymatrix@C=15pt{
(S_*,f)_{T_0} \ar@{..>}[dr] & \ar[l] (S_*,f)_{T_0\cap T_1}\ar@{..>}[d]  \ar[r] &  \ar@{..>}[dl] (S_*,f)_{T_1} \\
& (S_*,f)
}
\]
where the subscripts indicate evident pullback operations; for example $(S_*,f)_{T_0}$
is made from the portion of $S_*$ taken to $T_0$ by $f$. That diagram can be viewed
as a subcategory of $\eta/(S_*,f)$. It is a terminal subcategory in the sense that the inclusion
functor has a left adjoint. Hence the value of $\RR\eta_*(F_T|\mathscr V'_T)$ at the object $(S_*,f)$ of $\mathscr V_T$
can be identified with the homotopy pullback of
\[  F_T((S_*,f)_{T_0}) \rightarrow F_T((S_*,f)_{T_0\cap T}) \leftarrow F_T((S_*,f)_{T_1}) \]
which in turn can be identified with $F_T(S_*,f)=Z(\xi(S_*))$ by the sheaf property
of $Z$\,. This completes the proof.  \qed

\medskip
It is very easy to show (again) that the map~(\ref{eqn-denfinnrc}) is a weak equivalence using lemmas~\ref{lem-alpha} and~\ref{lem-xi}.
Indeed, lemma~\ref{lem-alpha} implies that the standard map
\[ \rmap(\xi^*N_d^{rc}P,\xi^*N_d^{rc}Q) \lra \rmap(\alpha^*\xi^*N_d^{rc}P,\alpha^*\xi^*N_d^{rc}Q) \]
is a weak equivalence and lemma~\ref{lem-xi} implies that the standard map
\[ \rmap(N_d^{rc}P,N_d^{rc}Q)\lra  \rmap(\xi^*N_d^{rc}P,\xi^*N_d^{rc}Q) \]
is a weak equivalence. See \cite[lem.A.1]{BoavidaWeissLong} and remember diagram~(\ref{eqn-decomp}).
This new proof carries over with a little mutatis mutandis to show that the
map~(\ref{eqn-denfinnrctrunc}) is also weak equivalence. In this way we have completed the proof of
theorem~\ref{thm-opersuppl}, since that had already been reduced to showing that~(\ref{eqn-denfinnrctrunc}) is a
weak equivalence.

\end{document}